\newif\ifSHOWEXTRA
\SHOWEXTRAtrue
\SHOWEXTRAfalse
\documentclass[12pt]{article}
\usepackage[utf8]{inputenc}
\usepackage[margin=1in]{geometry}
\usepackage{amsmath,amsthm,amssymb,amsfonts, pifont, bbm, xcolor, scrextend}
\definecolor{forestgreen}{RGB}{0,98,51}
\usepackage{appendix}
\usepackage{fancyhdr}
\usepackage{color,colortbl}
\newcounter{CommentCounter}

\definecolor{Gray}{gray}{0.9}
\definecolor{LightGreen}{rgb}{0,0.1,0.1}
\newcommand{\R}{\mathbb R}
\newcommand{\G}{\mathcal G}

\newcommand{\N}{\mathbb N}

\newcommand{\Z}{\mathbb Z}

\usepackage{tikz}
\usepackage{pgfplots}
\usepackage{eqnarray}
\usepackage{outlines}
\usepackage[mathscr]{euscript}
 \let\mathscr\relax
\usepackage[scr]{rsfso}
\usepackage{multicol}
\usepackage{enumerate}
\usepackage[colorlinks=true, pdfstartview=FitV, linkcolor=blue,
citecolor=blue, urlcolor=blue]{hyperref}

\newcommand{\totheleft}{-}
\newcommand{\totheright}{+}
\newcommand{\allowable}{\mathcal T}

\newtheorem{thm}{Theorem}[section]
\numberwithin{thm}{section}

\newtheorem{lem}[thm]{Lemma}

\newtheorem*{thm*}{Theorem}

\numberwithin{quest}{section}

\theoremstyle{remark}

\numberwithin{rem}{section}

\theoremstyle{definition}

\newtheorem{exmp}{Example}[section]

\newtheorem{claim}{Claim}
\numberwithin{claim}{exmp}

\pgfplotsset{compat=1.12}

\begin{document}
\author{
Daniel J. Slonim\footnote{Purdue University,
  Department of Mathematics,
150 N. University Street, West Lafayette, IN 47907, dslonim@purdue.edu, 0000-0002-9554-155X}
}
\title{On total weight exiting finite, strongly connected sets in shift-invariant weighted directed graphs on $\Z$}

\date{\today}

\maketitle

\begin{abstract}
    For a shift-invariant weighted directed graph with vertex set $\Z$, we examine the minimal weight $\kappa_0$ exiting a finite, strongly connected set of vertices. Although $\kappa_0$ is defined as an infimum, it has been shown that the infimum is always attained by an actual set of vertices.
    We show that for each underlying directed graph (prior to assignment of the weights), there is a formula for $\kappa_0$ as a minimum of finitely many integer combinations of the edge weights. We find this formula for several different directed graphs.
    Motivation for this problem comes from random walks in Dirichlet environments (equivalently, directed edge reinforced random walks), where the size of $\kappa_0$ has been shown to determine the strength of finite traps where the walk can get stuck for a long time.
    \\
    {\it Kewords:}
    directed graph,
    weighted directed graph,
    graph theory,
    Dirichlet environments,
    reinforced random walk,
\end{abstract}

\section{Introduction}

\begin{figure}
    \centering
    \includegraphics[width=6.5in]{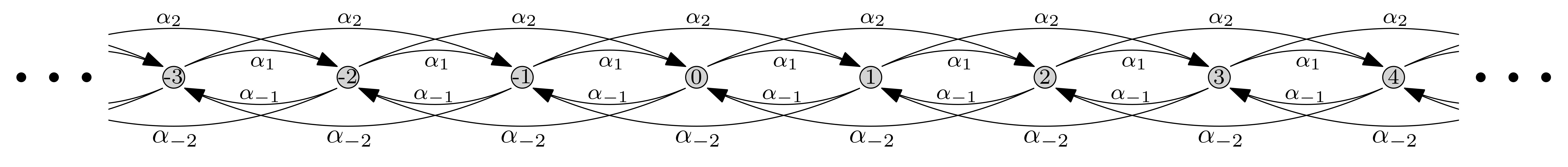}
    \caption{A portion of the graph $\G$ with $L=R=2$.}
    \label{fig:graphG}
\end{figure}

Let $\G=(\Z,E,w)$ be a weighted directed graph with vertex set $\Z$, edge set
$~E\subseteq \Z\times \Z$, and a weight function $w:E\to\R^{>0}$. If $e=(x,y)\in E$, we say that $e$ is an edge from $x$ to $y$, and we say the {\em head} of $e$ is $\overline{e}=y$ and the {\em tail} of $e$ is $\underline{e}=x$. We say a set $S\subset V$ is {\em strongly connected} if for all $x,y\in S$, there is a path from $x$ to $y$ in $\G$ using only vertices in $S$. Assume $\G$ satisfies the following assumptions.
\begin{enumerate}[(C1)]
    \item\hypertarget{cond:C1} $\G$ is shift-invariant.
    \item\hypertarget{cond:C2} $\G$ is strongly connected.
    \item\hypertarget{cond:C3} Each vertex of $\G$ has finite out-degree.
\end{enumerate}
These assumptions imply that there are positive integers $L$ and $R$ and nonnegative parameters $(\alpha_i)_{i=-L}^R$ with $\alpha_{-L},\alpha_{R}>0$, such that for $x,y\in\Z$, there is an edge in $\G$ from $x$
to $y$ if and only if $y\in\{x-L,\ldots,x+R\}$ and $\alpha_{y-x}>0$, and the weight of that edge, if it exists, is $\alpha_{y-x}$. 
An example with $L=R=2$ is represented in Figure \ref{fig:graphG} (here, and in other illustrations of graphs, we depict the case where $\alpha_0=0$, but our model does allow for $\alpha_0>0$). 

In this paper, we study the minimal weight exiting a finite, strongly connected subset of $\Z$. For a set $S\subseteq \Z$, define
\begin{equation}\label{betaS}
    \beta_S:=\sum_{\underline{e}\in S,~\overline{e}\notin S}w(e).
\end{equation}
This parameter $\beta_S$ is the sum of the weights of all edges exiting the set $S$. Define
\begin{equation}\label{eqn:kappa0def}
    \kappa_0:=\inf\{\beta_S:S\subset \Z\text{ finite, strongly connected}\}.
\end{equation}
Although $\kappa_0$ is defined as a minimum, the author showed in \cite{Slonim2021a} that it is actually a minimum: that is, it is attained as $\beta_S$ for some specific set $S$. That paper also gave an algorithm to find $\kappa_0$ given $L$, $R$, and the $\alpha_i$, $-L\leq i\leq R$. However, the algorithm requires examining a number of sets that grows exponentially with the inverse of the smallest positive $\alpha_i$\footnote{Even moderate examples can get unwieldy very quickly. For instance, say $L=5$, $R=4$, and $\alpha_{-5}=5$, $\alpha_3=0.1$, $\alpha_4=5$, and $\alpha_i=0$ for all other $i$. In this example, the $M$ calculated in \cite{Slonim2021a}, which is the largest diameter of a set one needs to check, is 4972. Further shrinking $\alpha_3$ increases this number still, so if one applies the algorithm as naively as possible, a moderate example can result in the need to examine on the order of $2^{5000}$ sets of vertices.}.

We show in this paper that for a given underlying directed graph, even before assignment of weights, there is a finite collection $\mathscr{S}$ of sets $S\subset\Z$ so that $\kappa_0$ is realized as $\beta_S$ for some $S\in\mathscr{S}$. Thus, there is a formula for $\kappa$ as a minimum of finitely many integer combinations of the $\alpha_i$. At present, we do not have a general algorithm to compute this formula for all possible underlying directed graphs, but we do it for several examples using {\em ad hoc} methods. These examples exhibit various important features that sets $S$ attaining $\beta_S=\kappa_0$ may have. 
Example \ref{kappa0example6} gives a particularly good demonstration of the difficulty that can attend such calculations.

\subsection{Motivation}

In \cite{Slonim2021a}, the author studied random walks in Dirichlet environments (RWDE) on $\Z$ with bounded jumps. These are random Markov chains on $\Z$ where transition probability vectors at each site are themselves chosen randomly, each being drawn according to a Dirichlet distribution. To the weighted directed graph $\G$, we can associate a RWDE model as follows. 

Recall the definition of the Dirichlet distribution: for a finite set $I$, take parameters ${\bf a}=(a_i)_{i\in I}$, with $a_i>0$ for all $i$. The Dirichlet distribution with these parameters is a probability distribution on the simplex $\Delta_I:=\{(x_i)_{i\in I}:\sum_{i\in I}x_i=1\}$ with density
\begin{equation*}D\left((x_i)_{i\in I}\right)=C({\bf a})\prod_{i\in I}x_i^{a_i-1},\end{equation*}
where $C({\bf a})$ is a normalizing constant. 

Let $\omega(x,y)$ denote the transition probability from $x$ to $y$. Let $\mathcal{N}\subseteq\{-L,\ldots,R\}$ be the set of $i$ for which $\alpha_i>0$. For each $x\in\Z$, let $(\omega(x,x+y))_{y\in\mathcal{N}}$ be drawn according to a Dirichlet distribution according to parameters $(\alpha_i)_{i\in\mathcal{N}}$. Doing this independently at each $x$ gives the RWDE model studied in \cite{Slonim2021a}. 

An interesting fact about this RWDE model is that due to a connection with Polya urns, they are equivalent to directed edge reinforced random walks (DERRWs) on $\G$. This is a walk whose transition rule is as follows: if the walk is at vertex $x$, then it jumps along an edge of $\G$ with probability proportional to the weight of that edge. However, each time the walk takes an edge, the weight of that edge is increased by 1, making it more likely that the walk will take the same edge the next time it is at site $x$. The law of a DERRW on $\G$ is the same as the law of the RWDE associated with the graph $\G$. This was first proven in \cite{Enriquez&Sabot2002} and \cite{Keane&Rolles2002}; see \cite{Sabot&Tournier2016} for useful discussion. 
 
The effect of reinforcement in DERRWs is stronger when the initial weights are small.
From the RWDE perspective, this corresponds to the fact that beta random variables\footnote{Dirichlet random vectors are a generalization of beta random variables, and thier individual components are distributed as beta random variables.} are more likely to take on values very close to 0 or 1 if they have small parameters. For this reason, smaller parameters mean the walk is more likely to get stuck in a finite set of vertices for a long time, iterating one or more cycles many times before escaping the set. Tournier \cite{Tournier2009} quantified this trapping effect in the following theorem.
\begin{thm*}[{\cite[Theorem 1]{Tournier2009}}]
Let $S$ be a finite, strongly connected subset of $\Z$ containing 0, and let $N_0^S({\bf X})$ denote the number of visits the walk ${\bf X}$ makes to 0 before leaving $S$. Then for any $r>0$,
\[
E[N_0^S]<\infty \quad \Leftrightarrow\quad\beta_S>1.
\]
\end{thm*}
This implies that if $\kappa_0\leq1$, then there is some set $S$ such that the expected number of visits to 0 before the walk leaves $S$ is infinite. The actual theorem from \cite{Tournier2009} is stronger, saying that for RWDE on a finite graph, a quantity known as the {\em quenched Green function} has finite moments up to (but not including) $\kappa_0$. 
On the infinite graph $\G$, the quenched Green function has finite moments up to (but not including) the minimum of $\kappa_0$ and another parameter $\kappa_1$, which is easily calculated from the $\alpha_i$, and the expected total number of visits to 0 is finite if and only if $\min(\kappa_0,|\kappa_1|)>1$. For more detials, see \cite{Slonim2021a}. Due to results from \cite{Slonim2022Ballistic}, the expected number of visits to 0 gives information about the {\em ballisticty} of the walk; if the expected number of visits is finite, then the walk is ballistic: that is, it has nonzero limiting speed. If the expected number of visits is finite, then the limiting speed is 0. 

\section{Theorems and Lemmas}

\begin{thm}\label{prop:GammaAttainedweakstrong}
Given $L$, $R$, and the $i$ for which $\alpha_i>0$, there is a collection $\mathscr{S}$ of finite, strongly connected sets $S\subseteq\Z$ such that $\kappa_0$ is attained as $\beta_S$ for some $S\in\mathscr{S}$.  
\end{thm}

The above theorem implies that for each underlying directed graph (prior to assignment of edge weights), there is a formula for $kappa_0$ as an elementary function of the $\alpha_i$. To prove it, we need the following lemma.

\begin{lem}\label{lem:SetTheory}
Let $\allowable\subseteq\N^{k}$ be a set of ordered $k$-tuples of positive integers. Let $\prec$ be the natural partial ordering on $\N^k$, $(n_1,\ldots,n_k)\prec(n_1',\ldots,n_k')$ if $n_i\leq n_i'$ for all $i$. Then there is a finite subset $\allowable^*\subseteq\allowable$ such that for all ${\bf x}\in\allowable$, there is an ${\bf x}^*\in\allowable^*$ such that ${\bf x}^*\prec{\bf x}$.
\end{lem}

\begin{proof}
We prove this by induction on $k$. The base case $k=1$ is trivial; a 1-tuple is simply a positive integer, and we can let $\allowable^*:=\{\min\allowable\}$.
Now suppose the result is true for all subsets of $\N^k$, and let $\allowable\subseteq \N^{k+1}$. Now let $\underline{\allowable}:=\{(n_1,\ldots,n_k):(n_1,\ldots,n_k,n_{k+1}\in\allowable\text{ for some }n_{k+1}\in\N\}$ be the projection of $\allowable$ onto $\N^k$, and for $n\in\N$ define $\underline{\allowable}(n):=\{(n_1,\ldots,n_k):(n_1,\ldots,n_k,n)\in\allowable\}$. Thus, $\underline{\allowable}=\bigcup_{n=1}^{\infty}\underline{\allowable}(n)$.  
Now by the inductive hypothesis, there is a finite set $\underline{\allowable}^*\subseteq \underline{\allowable}$ such that every element of $\underline{\allowable}$ is greater than some element of $\underline{\allowable}^*$. Since $\underline{\allowable}^*$ is finite, $\underline{\allowable}^*\subseteq\bigcup_{n=1}^{N}\underline{\allowable}(n)$ for some $N$. Applying the inductive hypothesis to each $\underline{\allowable}(n)$ gives us sets $\underline{\allowable}^*(n)$ such that for all ${\bf x}\in\underline{\allowable}(n)$, there is a ${\bf x}^*\in\underline{\allowable}^*(n)$ such that ${\bf x}^*\prec{\bf x}$.

Now define $\allowable^*:=\bigcup_{n=1}^N\{(n_1,\ldots,n_k,n):(n_1,\ldots,n_k)\in\underline{\allowable}^*(n)\}$. It is easy to see that this is a finite subset of $\allowable$. Now suppose $(n_1,\ldots,n_k,n_{k+1})\in\allowable$. 
Suppose $n_{k+1}< N$. Then $(n_1,\ldots,n_k)\in\underline{\allowable}(n_{k+1})$, so there exists $(n_1^*,\ldots,n_k^*)\in\underline{\allowable}^*(n_{k+1})$ with $n_i^*\leq n_i$ for $1\leq i\leq k$. Since $(n_1^*,\ldots,n_k^*,n_{k+1})\in\allowable^*$, we are done. 

On the other hand, suppose $n_{k+1}\geq N$. Then since $(n_1,\ldots,n_k)\in\underline{\allowable}$, there exists $(n_1^*,\ldots,n_k^*)\in\underline{\allowable}^*$ such that $n_i^*\leq n_i$ for all $i=1,\ldots,k$. Now $(n_1^*,\ldots,n_k^*)\in\underline{\allowable}(n)$ for some $n\leq N\leq n_{k+1}$. Hence there exists $(n_1^{**},\ldots,n_k^{**})\in\underline{\allowable}^*(n)$ such that $n_i^{**}\leq n_i^*$ for all $1\leq i\leq k$. Thus $(n_1^{**},\ldots,n_k^{**},n)\in\allowable^*$ with $n_i^{**}\leq n_i$ for $1\leq i\leq k$, and $n\leq n_{k+1}$. \end{proof}

We are now able to prove Theorem \ref{prop:GammaAttainedweakstrong}.

\begin{proof}[Proof of Theorem \ref{prop:GammaAttainedweakstrong}]
For any finite set $S\subset \Z$, $\beta_S$ is a sum of weights of edges exiting $\beta$. The weight of each edge is $\alpha_i$ for some $i$, and each $\alpha_i$ must be included at least once, as the weight of an edge exiting either the rightmost or leftmost point of $S$. Thus,
\begin{equation*}
    \beta_S=\sum_{i=-L}^Rx_i\alpha_i
\end{equation*}
where $x_i=x_i(S):=\#\{z\in S:z+i\notin 
S\}\geq1$.

Now let $\allowable\subset\N^{R+L+1}$ be the set of ordered tuples $(y_{-L},\ldots,y_R)$ such that there is some finite, strongly connected set $S$ with $x_i(S)=y_i$ for all $-L\leq i\leq R$. Thus,
\begin{equation}\label{eqn:247}
    \kappa_0=\inf\left\{\sum_{i=-L}^Ry_i\alpha_i:(y_{-L},\ldots,y_R)\in\allowable\right\}.
\end{equation}
 Applying Lemma \ref{lem:SetTheory}, we get a finite set $\allowable^*\subseteq\allowable$ such that for any $S$, there is a $(y_{-L},\ldots,y_R)\in\allowable^*$ with $y_i\leq x_i(S)$ for all $-L\leq i\leq R$. This $(y_{-L},\ldots,y_R)$ is $(x_{-L}(S'),\ldots,x_R(S'))$ for some finite, strongly connected $S'$, so that $x_i(S')\leq x_i(S)$ for all $i$ and therefore $\beta_{S'}\leq S$. It follows that $\kappa_0=\beta_{S'}$ for some such $S'$. Choosing one such $S'$ for each element of $\allowable^*$ gives the finite set $\mathscr{S}$ described in the statement of the theorem.
\end{proof}

We remark that our proof allows us to rewrite \eqref{eqn:247} as 
\begin{equation}\label{kappa0elementaryfunction}
    \kappa_0=\min\left\{\sum_{i=-L}^Ry_i\alpha_i:(y_{-L},\ldots,y_R)\in\allowable^*\right\},
\end{equation}
giving $\kappa_0$ as an elementary function (a minimum of finitely many integer combinations) of the $\alpha_i$. However, because we do not have a constructive way to find the set $\allowable^*$, we do not have a general algorithm to find this formula for any given underlying directed graph. We are able to find the formula in several examples, which we do in Section \ref{sec:Examples}.

We  recall the following lemma, which is Claim 5.1.2 from the proof of Theorem 5.1 in \cite{Slonim2021a}. We include the proof to make this article more self-contained.

\begin{lem}\label{claim:GammaClaim3}
 Let $S\subset\Z$ be a finite, strongly connected set of vertices. If $x$ is a vertex to the left or to the right of $S$, then $\beta_{S\cup\{x\}}\geq\beta_{S}$. 
\end{lem}

\begin{proof}
Let $c^{\totheright}:=\sum_{i=1}^R\alpha_i$ and $c^{\totheleft}:=\sum_{i=-L}^{-1}\alpha_i$. 
The quantity $\beta_S$ is the sum of all weights from vertices in $S$ to vertices not in $S$. The quantity $\beta_{S\cup\{x\}}$ counts all same weights, except for weights of edges from $S$ to $x$, and it also counts weights of edges from $x$ to vertices not in $S\cup\{x\}$. 
If $x$ is to the right of $S$, then the total weight of edges from $S$ to $x$ cannot be more than $c^{\totheright}$, because $c^{\totheright}$ is the total weight into $x$ from all vertices to the left of $x$. On the other hand, $c^{\totheright}$ is also the total weight from $x$ to all vertices to the right of $x$, which are necessarily not in $S\cup\{x\}$. Thus, the additional weight from $x$ to the right at least makes up for any weight into $x$ from $S$. This proves the claim in the case that $x$ is to the right of $A$, and a similar argument proves the symmetric case.
\end{proof}

Using the notation $c^{\totheright}$ and $c^{\totheleft}$ introduced in the above proof, we note that 
\begin{equation}\label{eqn:boundsforkappa0}
    \kappa_0\geq c^{\totheright}+c^{\totheleft}.
\end{equation}
This is because any strongly connected set will have weight at least $c^{\totheright}$ exiting from the rightmost point and weight at least $c^{\totheleft}$ exiting from the leftmost point. Therefore, every strongly connected set $S$ has $\beta_S\geq c^{\totheright}+c^{\totheleft}$, and taking the infimum preserves the inequality.

\section{Examples}\label{sec:Examples}

We now give examples where we can find the formula for $\kappa_0$. Recall that $\kappa_0:=\inf\{\beta_S:S\subset \Z\text{ finite, strongly connected}\}$, where $\beta_S$ is the sum of edge weights leaving the set $S$. By shift invariance of the graph $\G$, it suffices to consider sets $S$ whose leftmost point is 0. 

\begin{exmp}\label{kappa0example1}
$L=R=1$.

In this case, $\kappa_0=\alpha_1+\alpha_{-1}$. This is because the only strongly connected sets are intervals, which all have exit weight.
\end{exmp}

\begin{exmp}\label{kappa0example7}
$\alpha_0>0$. 
 
In this case, $\{0\}$ is already a strongly connected set, so \eqref{eqn:boundsforkappa0} gives $\kappa_0=\beta_{\{0\}}=c^{\totheright}+c^{\totheleft}$. 
\end{exmp}

\begin{exmp}\label{kappa0example2}
$L=2$, $R=3$, $\alpha_i=0$ for $i=-1,\ldots,2$. 

In this case, we also have $\kappa_0=d^{\totheright}+d^{\totheleft}$. Let $S$ be a strongly connected finite set of vertices with left endpoint 0. Then $S$ contains 3, since there must be a vertex reachable from 0 in one step, and by assumption there are no vertices to the left of 0. Also, $S$ contains 2, since 0 must be reachable in one step from a vertex in $S$. Likewise, 2 must then also be reachable, and since $-1\notin S$, $S$ must contain 4 as well. Now since a vertex must be reachable from 3, $S$ must contain either 1 or 6. Suppose $S$ contains 1. Then $S$ contains $[0,4]$, which has exit weight $d^{\totheright}+d^{\totheleft}$, and by Lemma \ref{claim:GammaClaim3}, $\beta_S\geq \beta_{[0,4]}=d^{\totheright}+d^{\totheleft}$. On the other hand, suppose $S$ does not contain 1. Then it contains 6. If $S$ also contains 5, then $S$ contains the interval $[2,6]$, which is shift-equivalent to $[0,4]$. If $S$ contains neither 1 nor 5, then it contains exactly the set $\{0,2,3,4,6\}$ and possibly vertices to the right and/or left of this set, so by Lemma \ref{claim:GammaClaim3}, $\beta_S\geq\beta_{\{0,2,3,4,6\}}$. One can easily check that in this case, $\beta_{\{0,2,3,4,6\}}=d^{\totheright}+d^{\totheleft}$. 

We have calculated $\kappa_0$ without even showing that either $[0,4]$ or $\{0,2,3,4,6\}$ is strongly connected, but in fact they both are. Consider the path $0\to3\to1\to4\to2\to0$ in $[0,4]$ and the path $0\to3\to6\to4\to2\to0$ in $\{0,2,3,4,6\}$. Thus, even when $\kappa_0=d^{\totheright}+d^{\totheleft}$, a minimizing set $S$ for $\beta_S$ need not be an interval (although a large enough interval will always be a minimizing set). 
\end{exmp}

\begin{exmp}\label{kappa0example3}
$L=1$, $R\geq2$, $\alpha_0=0$ $\alpha_1>0$, $\alpha_i=0$ for $i=2,\ldots, R-1$.

In this case, we show that $\kappa_0=2\alpha_R+\alpha_1+\alpha_{-1}$; thus if $R>2$, then $\kappa_0<d^{\totheright}+d^{\totheleft}$. Let $S$ be a strongly connected set with left endpoint 0. Since 0 must be reachable from another point in $S$, we have $1\in S$. Now by Lemma \ref{claim:GammaClaim3}, this implies $\beta_S\geq\beta_{\{0,1\}}=2\alpha_R+\alpha_1+\alpha_{-1}$. The set $\{0,1\}$ is strongly connected, and hence $\kappa_0=\beta_{\{0,1\}}=2\alpha_R+\alpha_1+\alpha_{-1}$. 
\end{exmp}




\begin{exmp}\label{kappa0example4}
$L=R=2$, $\alpha_{-1},\alpha_1>0$, $\alpha_0=0$.

There are two possibilities. Let $S$ be a strongly connected set with leftmost point 0. If $1\in S$, then by Lemma \ref{claim:GammaClaim3}, $\beta_S\geq\beta_{\{0,1\}}=d^{\totheright}+d^{\totheleft}$. On the other hand, if $1\notin S$, then $2\in S$, and by Lemma \ref{claim:GammaClaim3}, $\beta_S\geq\beta_{\{0,2\}}=\alpha_{-2}+2\alpha_{-1}+2\alpha_1+\alpha_2$. Since both $\{0,1\}$ and $\{0,2\}$ are strongly connected, $\kappa_0$ may be either $d^{\totheright}+d^{\totheleft}$ or $\alpha_{-2}+2\alpha_{-1}+2\alpha_1+\alpha_2$, depending on whether $\alpha_{-1}+\alpha_1$ or $\alpha_{-2}+\alpha_2$ is smaller. That is, 
\begin{align*}
        \kappa_0&=\min(2\alpha_{-2}+\alpha_{-1}+\alpha_1+\alpha_2,\alpha_{-2}+2\alpha_{-1}+2\alpha_1+\alpha_2)
        \\
        &=\alpha_{-2}+\alpha_{-1}+\alpha_1+\alpha_2+\min(\alpha_{-1}+\alpha_1,\alpha_{-2}+\alpha_2)
\\
&=\min(\beta_{\{0,1\}},\beta_{\{0,2\}})
. 
\end{align*}
\end{exmp}

\begin{exmp}\label{kappa0example5}
$L=6$, $R=3$, $\alpha_2>0$, $\alpha_i=0$ for $i=-5,\ldots,1$.  

If $S$ is a finite, strongly connected set with 0 the leftmost vertex, then $6\in S$, since $0$ must be reachable from the right. We consider possible sets $S\cap[0,6]$. There are 32 subsets of $[0,6]$ that contain 0 and 6; however, $S$ must contain either 2 or 3, since there must be edges from 0 to other sets in $S$ and nothing to the left of 0 is allowed. Similarly, if $S$ contains 1, then it must contain either 3 or 4, and if $S$ contains 2, then it must contain either 4 or 5. This eliminates 12 of the 32 possibilities, leaving 20 possibilities for $S\cap[0,6] $. Of these, we first consider two candidates, $\{0,3,6\}$ and $\{0,2,4,6\}$. Both of these are strongly connected, and so $\beta_{\{0,3,6\}}=2\alpha_{-6}+3\alpha_2+\alpha_3$ and $\beta_{\{0,2,4,6\}}=3\alpha_{-6}+\alpha_2+4\alpha_3$ both provide upper bounds for $\kappa_0$. Depending on the values of the $\alpha_i$, either can be lower than the other. The set $\{0,2,3,4,6\}$ is also strongly connected, but has $\beta_{\{0,2,3,4,6\}}=4\alpha_{-6}+2\alpha_2+3\alpha_3$. Thus, if $\alpha_2\geq\alpha_3$, then $\beta_{\{0,2,4,6\}}<\beta_{\{0,2,3,4,6\}}$, and if $\alpha_2\leq\alpha_3$, then $\beta_{\{0,3,6\}}<\beta_{\{0,2,3,4,6\}}$. One can simply check that other 17 of the possible sets $D=S\cap[0,6]$ either have $\beta_D>\beta_{\{0,3,6\}}$ for all possible values of the $\alpha_i$, 
$\beta_D>\beta_{\{0,2,4,6\}}$ for all possible values of the $\alpha_i$,
or
$\beta_D>\beta_{\{0,2,3,4,6\}}$ for all possible values of the $\alpha_i$. By Lemma \ref{claim:GammaClaim3}, this implies that $\beta_S\geq\min(\beta_{\{0,3,6\}},\beta_{\{0,2,4,6\}})$. Therefore, 
\begin{align*}
    \kappa_0
    &=\min(2\alpha_{-6}+3\alpha_2+\alpha_3,3\alpha_{-6}+\alpha_2+4\alpha_3)
    \\
    &=\min(\beta_{\{0,3,6\}},\beta_{\{0,2,4,6\}}).
\end{align*}
\end{exmp}

In all five of the above examples, there is always a set $S$ minimizing $\beta_S$ that represents a single, simple loop. The exit time from $S$ is the first time the walk stops repeating this loop. Thus, if $\kappa_0\leq1$, then there is a single loop that the walk is expected to repeat infinitely many times before deviating from it. 

In the nearest-neighbor case, treated in Example \ref{kappa0example1}, $\kappa_0\leq1$ means the walk is expected to repeat the loop $0\to1\to0$ infinitely many times before ever taking a different step (and, likewise, the walk is expected to repeat the loop $0\to-1\to0$ infinitely many times before ever stepping to $1$). This does not mean the only finite traps are sets of the form $\{x,x+1\}$. For example, it is also the case that $\beta_{[0,5]}\leq1$, so that the walk is expected to spend an infinite amount of time in $[0,5]$ before leaving it, regardless of the precise path (and even if transition probabilities at sites 1,2,3, and 4 are conditioned to be moderate). But there are no finite traps ``worse" (in the sense of finite moments of quenched expected exit time) than the set $\{0,1\}$. 

In fact, for nearest-neighbor RWDE on $\Z^d$, pairs of adjacent vertices are always the worst finite traps, and if $\kappa_0\leq1$, then the walk is expected to bounce back and forth between 0 and one other vertex infinitely many times before doing anything else \cite{Sabot&Tournier2016}.

Our other examples so far match this trend in a sense; although the worst traps are not necessarily pairs of vertices, the worst traps are loops, and $\kappa_0\leq1$ means there is a loop that the walk is expected to iterate infinitely many times before doing anything else. 
\begin{itemize}
    \item In Example \ref{kappa0example7}, one such loop is $0\to0$. 
    \item In Example \ref{kappa0example2}, one such loop is $0\to3\to6\to4\to2\to0$. 
    \item In Example \ref{kappa0example3}, one such loop is $0\to1\to0$. 
    \item In Example \ref{kappa0example4}, one such loop is $0\to1\to0$ (if $\beta_{\{0,1\}}\leq1$) or $0\to2\to0$ (if $\beta_{\{0,2\}}\leq1$).
    \item In Example \ref{kappa0example5}, one such loop is $0\to3\to6\to0$ (if $\beta_{\{0,3,6\}}<1$) or $0\to2\to4\to6$ (if $\beta_{\{0,2,4,6\}}<1$). 
\end{itemize}

Our next example shows that unlike in the nearest-neighbor case on $\Z^d$, there are parameters where the strongest finite traps never represent just one loop. In particular, one can find cases where $\kappa_0\leq1$, so there are finite traps in which the walk is expected to be stuck for an infinite amount of time, but there is no single loop that the walk is expected to iterate infinitely many times before deviating from it.

\begin{exmp}\label{kappa0example9}
$L=R=2$, $\alpha_{-1}=\alpha_0=0$, $\alpha_1>0$.
\end{exmp}

A finite, strongly connected set with 0 as its leftmost point will necessarily contain 2, since 0 must be reachable from the right. Thus, by Lemma \ref{claim:GammaClaim3}, for any finite, strongly connected $S$, $\beta_S\geq\min(\beta_{\{0,1,2\}},\beta_{\{0,2\}})$. Now $\{0,2\}$ and $\{0,1,2\}$ are already strongly connected, so $\kappa_0=\min(\beta_{\{0,1,2\}},\beta_{\{0,2\}})$. The minimum may be achieved on either set, depending on the $\alpha_i$. 

\begin{figure}[ht]
    \centering
    \includegraphics[width=6.5in]{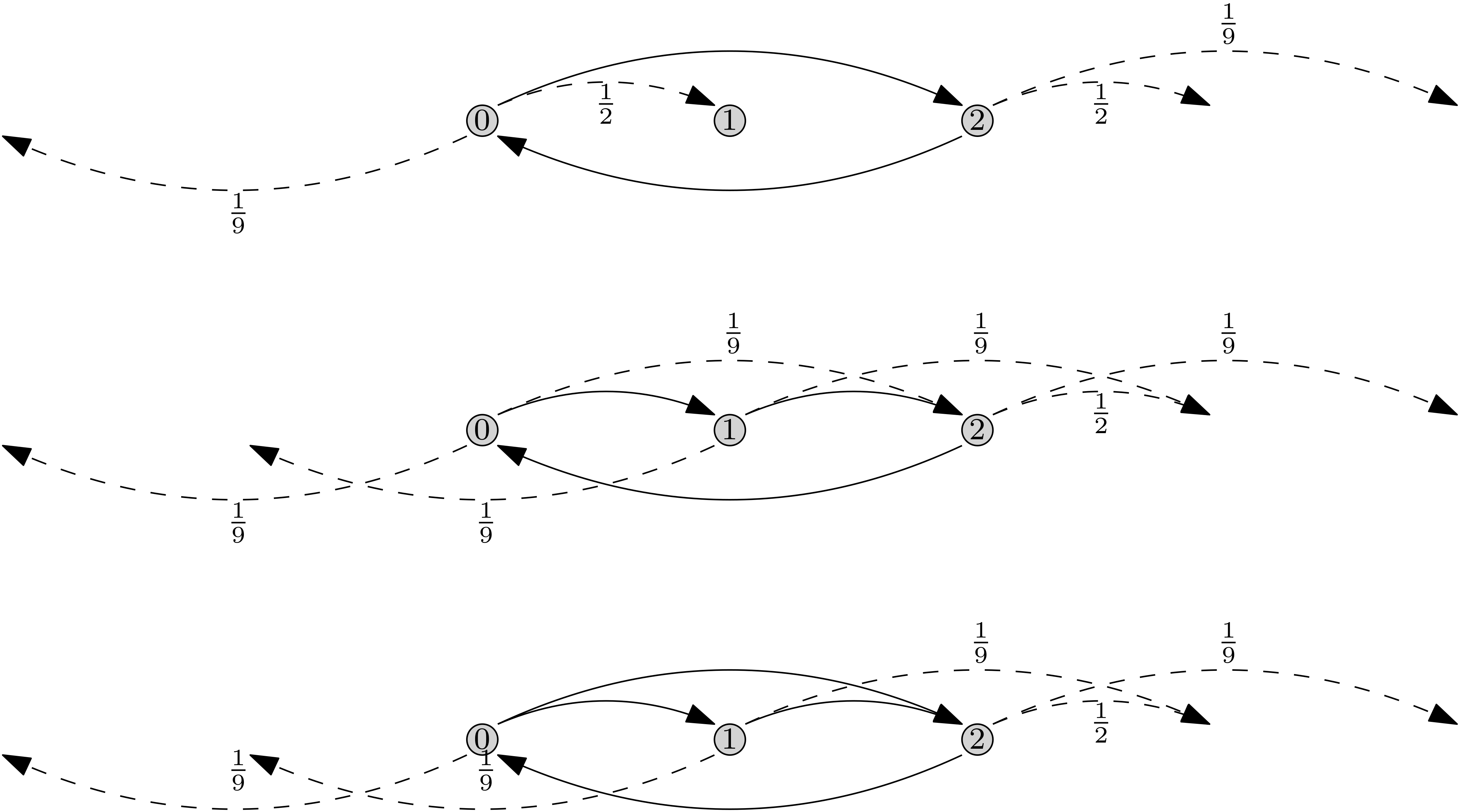}
    \caption{The top shows the weights exiting the loop $0\to2\to0$. The middle shows the weights exiting the loop $0\to1\to2\to0$. The bottom shows the weights exiting the union of these two loops, or the set $\{0,1,2\}$.}
    \label{fig:kappa0figure}
\end{figure}

We now examine a case where $\kappa_0\leq1$, but there are no loops that the walk is expected to iterate infinitely many times before doing anything else. Suppose $\alpha_{-2}=\alpha_{2}=\frac19$, and $\alpha_1=\frac12$. Then $\beta_{\{0,2\}}=\frac{11}{9}>1$, and $\beta_{\{0,1,2\}}=\frac{17}{18}<1$. Thus, $\kappa_0=\frac{17}{18}$, and a walk started from 0 is expected to spend an infinite amount of time in $\{0,1,2\}$ before exiting. However, because $\beta_{\{0,2\}}=\frac{11}{9}>1$, the expected exit time from $\{0,2\}$ is finite, so the walk is not expected to iterate the loop $0\to2\to0$ infinitely many times before deviating from it. Moreover, the walk is not expected to iterate the loop $0\to1\to2\to0$ infinitely many times before deviating it, but to see this, we must use the original formulation of Tournier's lemma from \cite{Tournier2009}. The formulation there is in terms of sets of edges rather than vertices. The edges that are not in the loop $0\to1\to2$ but have tails in the vertex set touched by this loop have weights that add up to $\frac{19}{18}>1$. Hence \cite[Theorem 1]{Tournier2009} implies that the expected time to deviate from this set of edges (and thus from the loop $0\to1\to2\to0$) is finite. Nevertheless, the weight exiting $\{0,1,2\}$ is $\frac{17}{18}<\frac12$, so the walk is expected to stick to the vertex set $\{0,1,2\}$, and thus to the pair of loops $0\to2\to0$ and $0\to1\to2\to0$, for an infinite amount of time before doing anything else. See Figure \ref{fig:kappa0figure}.

Note that in this example, $\beta_{\{0,1,2\}}<\beta_{\{0,2\}}$. This shows that in Lemma \ref{claim:GammaClaim3}, the assumption that $x$ is to the right or left of $x$ is really needed.
 
Our next example presents a similar phenomenon: the walk is not expected to get stuck in any one loop for an infinite amount of time, but the walk is expected to spend an infinite amount of time in a set of vertices. In the previous example, the vertex set $S$ minimizing $\beta_S$ can have all of its vertices hit by one loop, the loop $0\to1\to2\to0$, but this loop alone does not have a trapping effect as strong as the whole set $S$. In the next example, there is no single loop that can hit all the vertices in the minimizing $S$, so our formulation of Tournier's lemma in terms of vertices only is enough to see that there is no single loop the walk is expected to traverse infinitely many times before straying from it. The next example also presents a calculation of $\kappa_0$ for a situation where it is less straightforward than the others we've examined.

\begin{exmp}\label{kappa0example6}
$L=16$, $R=5$, $\alpha_{-16},\alpha_2,\alpha_5>0$, all other $\alpha_i=0$. 

In this case, there are four possible values for $\kappa_0$, three of which can be attained by sets of vertices representing single loops, but one of which cannot. We will show that $\kappa_0$ is attained by one of the following four sets:

\begin{itemize}
    \item ${S_1}=\{0,2,4,6,8,10,12,14,16\}$. This set represents a loop that steps up by 2s from 0 to 16 and then jumps back to 0. $\beta_{S_1}=8\alpha_{-16}+   \alpha_2+9\alpha_5$.
    \item ${S_2}=\{0,5,10,15,16,20,25,30,32\}$.  The set ${S_2}$ represents a loop that steps up by 5s from 0 to 30, then steps to 32 and jumps back to 16 and then to 0. This is one of 28 loops that all step up by 5 six times, up by 2 once, and down by 16 twice, having vertex set $S$ with leftmost point 0. All such loops have the same associated $\beta_{S}=\beta_{S_2}=7\alpha_{-16}+   8\alpha_2+     3\alpha_5$. 
    \item $S_3=\{0,5,10,12,14,16\}$. The set ${S_3}$ represents a loop that steps up by 5s from 0 to 10, then by 2s from 0 to 16, then jumps back to 0. This is one of 10 loops that all have vertex set $S\subset[0,16]$, all of which have $\beta_{S}=\beta_{S_3}=    5\alpha_{-16}+   3\alpha_2+     4\alpha_5$.
    \item ${S_4}=\{0,2,4,5,6,7,8,9,10,11,12,14,16\}$. This set does {\em not} represent one single loop; in fact, it represents all 10 loops that stay within $[0,16]$.  $\beta_{S_4}=    12\alpha_{-16}+   2\alpha_2+     5\alpha_5$.
\end{itemize}
One can check that any of $\beta_{S_1}$, $\beta_{S_2}$, $\beta_{S_3}$, and $\beta_{S_4}$ can be the smallest, depending on the $\alpha_i$. We will show
\begin{equation*}
    \kappa_0=\min(\beta_{S_1},\beta_{S_2},\beta_{S_3},\beta_{S_4}). 
\end{equation*}

To confirm that these are the possible values for $\kappa_0$, let $S$ be a finite, strongly connected set with leftmost point 0, and we will show that $\beta_S$ is at least as large as one of these values. First, we note
\begin{equation*}
    \beta_S=x_{-16}\alpha_{-16}+x_2\alpha_2+x_5\alpha_5,
\end{equation*}
where $x_i=x_i(S):=\#\{z\in S:z+i\notin S\}\geq1$. Since 0 must be reachable in one step from a vertex to its right, $16\in S$.

\begin{claim}\label{claim:16proofclaim1}
$x_{-16}\geq 5$.
\end{claim}

$S$ must contain $0$ and $16$, so there must be a path $\sigma$ from $0$ to $[16,\infty)$ that does not leave $S$. Since the only step to the left is down 16, all steps in this path must be to the right (and must therefore be of length 2 or 5). If this path includes two or more steps of length 2, then $S\cap[0,15]$ must have at least 5 elements. But for each $z\in S\cap[0,15]$, $z-16\notin S$, so each element of $S\cap[0,15]$ contributes 1 to $x_{-16}$. Hence $x_{-16}\geq5$. On the other hand, if the path $\sigma$ includes no steps or one step of length 2, then the path $\sigma$ includes four vertices in $S\cap[0,15]$, does not include 1 or 4, and lands on either 17 or 20. Since 1 and 4 are not in $\sigma$, either 17 or 20 contributes 1 to $x_{-16}$, in which case we then have $x_{-16}\geq5$, or else 1 or 4 is in $S$, in addition to the four vertices from $\sigma$ that are in $[0,15]$, so that $|S\cap[0,15]|\geq5$, and so $x_{-16}\geq5$. This proves our claim. 

\begin{claim}\label{claim:16proofclaim2}
If $x_2=1$, then $x_{-16}\geq8$ and $x_5\geq9$.
\end{claim}

To see this, note that if $x_2=1$, then $S$ includes only even vertices; otherwise, the rightmost odd vertex and the rightmost even vertex would each contribute 1 to $x_2$, giving $x_2\geq2$. Moreover, since $S$ contains $0$ and $16$, it must contain every even vertex between, in order to prevent any vertex other than the rightmost from contributing to $x_2$. The 8 even vertices $z=0$ through $z=14$ each have $z-16\notin S$, so $x_{-16}\geq8$, and the 9 even vertices $z=0$ through $z=16$ each have $z+5\notin S$, so $x_5\geq9$. 

\begin{claim}\label{claim:16proofclaim3}
$x_5\geq3$. 
\end{claim}

To see this, note that the rightmost vertex from each equivalence class $\pmod 5$ will contribute 1 to $x_5$. We already have $0,16\in S$, so the equivalence classes $0$ and $1$ are represented. But the equivalence class 1 is only reachable from equivalence class 2 (via a downward step of length 16) and from equivalence class 4 (via an upward step of length 2). Hence $S$ must contain an element from one of the equivalence classes 2 or 4 $\pmod 5$, and therefore at least three equivalence classes are represented, so $x_5\geq3$. 

\begin{claim}\label{claim:16proofclaim4}
If $x_5=3$, then $x_{-16}\geq7$, and $x_2\geq8$. 
\end{claim}

We first note that since each equivalence class contributes only 1 to $x_5$, all elements in each equivalence class must form an unbroken arithmetic progression from the lowest to the highest. That is, letting $z_i^{\text{least}}$ and $z_i^{\text{greatest}}$ be, respectively, the least and greatest $z$ such that $z\equiv i \pmod 5$ and $z\in S$, we have $\{z_i^{\text{least}},z_i^{\text{least}}+5,z_i^{\text{least}}+10,\ldots,z_i^{\text{greatest}}\}\subset S$. We now examine two separate cases.

{\em Case 1: $S$ contains elements from equivalence classes 0,1, and 2 $\pmod5$.}

Since $S$ contains no elements from equivalence class 4, equivalence $z_1^{\text{least}}$ can only be reached from equivalence class 2, which occurs via a leftward step of length $16$. Thus $z_1^{\text{least}}+16\in S$.
On the other hand, since $S$ contains no elements from equivalence class 3, equivalence class 2 can only be reached from equivalence class 0, via a rightward step of length 2. Thus $z_2^{\text{least}}-2\in S$.

Therefore, the path
\begin{equation*}
    0\to5\to10\to\cdots\to (z_2^{\text{least}}-2)\to z_2^{\text{least}}
    \to\cdots\to (z_1^{\text{least}}+16)\to z_1^{\text{least}}\to\cdots\to 16\to0
\end{equation*}
is in $S$. All steps marked out by ellipses are upward steps of length 5. It follows that this is a path of length 9, since any other number of steps would result in ending at a point other than 0. 

All the vertices in equivalence class 0 contribute 1 to $x_{-16}$, since stepping down by 16 would reach a vertex in equivalence class 4, which cannot be in $S$. Vertices from the path that are in equivalence class 2, other than $z_1^{\text{least}}+16$, are less than $z_1^{\text{least}}+16$, and so stepping down by 16 reaches a vertex that is in equivalence class 1 but not in $S$. And all vertices from the path that are in equivalence class 1, other than 16, are less than 16, so stepping down by 16 reaches a vertex not in $S$. Thus all but two of the vertices from the path shown will contribute 1 to $x_{-16}$, and therefore $x_{-16}\geq7$.

Now all the vertices in equivalence classes 1 or 2 contribute 1 to $x_2$, since stepping to the right by 2 reaches a vertex in equivalence class 3 or 4. And vertices in equivalence class 0 that are less than $z_2^{\text{least}}-2$ also contribute to $x_2$, since stepping to the right by 2 reaches a vertex in equivalence class 2 but less than $z_2^{\text{least}}$. Thus, all but one of the vertices shown in this path contribute to $x_2$, so $x_2\geq8$.

{\em Case 2: $S$ contains elements from equivalence classes 0,1, and 4 $\pmod5$.}

By a similar argument to that given in Case 1, the path
\begin{equation*}
    0\to5\to10\to\cdots\to (z_4^{\text{least}}+16)\to z_4^{\text{least}}
    \to\cdots\to (z_1^{\text{least}}-2)\to z_1^{\text{least}}\to\cdots\to 16\to0
\end{equation*}
is in $S$. All steps marked out by ellipses are upward steps of length 5. It follows that this is a path of length 9, since any other number of steps would result in ending at a point other than 0. Now, for $z=0,5,10,z_4^{\text{least}}+11$, we have $z-16\equiv 4\pmod 5$, but $z-16<z_4^{\text{least}}$, so $z-16\notin S$ and $z$ contributes 1 to $x_{-16}$. Moreover, $z_4^{\text{least}}$ and every subsequent vertex are all less than 16 (except, of course, for 16 itself), so they all contribute 1 to $x_{-16}$. Thus, $x_{-16}\geq7$.

Moreover, all vertices in equivalence class 0 or 1 contribute 1 to $x_2$, since $S$ has no vertices in equivalence class 2 or 4. And all but one of the vertices $z$ in equivalence class 4 are strictly less than $z_1^{\text{least}}-2$, so that $z+2\notin S$. Hence all but one of the vertices in the loop contribute to $x_2$, so $x_2\geq8$.

\begin{claim}\label{claim:16proofclaim5}
If $x_2=2$, then $x_{5}\geq5$. 
\end{claim}
If $x_2=2$, then $S$ contains even and odd elements (because the only even upward jumps are of length 2, a strongly connected $S$ with only even elements would have $x_2=1$). It therefore must contain every even number, from its least even number to its greatest even number. In particular, it must contain $S_1=[0,16]$. This is enough to include at least one representative from every equivalence class $\pmod 5$. The greatest element of $S$ in each of these equivalence classes contributes 1 to $x_5$, so $x_5\geq5$. 

\begin{claim}\label{claim:16proofclaim7}
If $x_2=2$, then $x_5+x_{-16}\geq 17$ and $x_{-16}\geq9$. 
\end{claim}

We have already established that if $x_2=2$, then $S$ contains $S_1=[0,16]$ and at least one odd number.
Now $S_1$ has $x_{-16}=8$ and $x_5=9$. The odd number will also contribute 1 to $x_{-16}$, giving the bound $x_{-16}\geq9$. 
The set $S_1$ includes 8, which is in equivalence class 3 $\pmod 5$, and two elements of each of the equivalence classes 0,1,2, and 4. Each of the equivalence classes must contribute at least 1 to $x_5$, and for any of the classes 0,1,2, or 4 to avoid contributing 2, the odd number in between the two even numbers from that equivalence class must be contained in $S$. This saves 1 from $x_5$ but adds 1 to $x_{-16}$, thus keeping $x_5+x_{-16}\geq17$.

\begin{claim}\label{claim:16proofclaim8}
$\beta_S\geq\min(\beta_{S_1},\beta_{S_2},\beta_{S_3},\beta_{S_4})$. 
\end{claim}

We know $\beta_S$ must have $x_5\geq3$ by Claim \ref{claim:16proofclaim3}. By Claim \ref{claim:16proofclaim4}, if $x_5=3$, then $\beta_S\geq7\alpha_{-16}+8\alpha_2+3\alpha_5=\beta_{S_2}$. Now suppose $x_5\geq 4$. If $x_2=1$, then $\beta_S\geq 8\alpha_{-16}+\alpha_2+9\alpha_5=\beta_{S_1}$ by Claim \ref{claim:16proofclaim2}. Now consider the case $x_2=2$. Then $x_5+x_{-16}\geq17$ by Claim \ref{claim:16proofclaim7}. If $\alpha_5>\alpha_{-16}$, then since $x_5\geq5$ by Claim \ref{claim:16proofclaim5}, we have $\beta_S\geq 12\alpha_{-16}+2\alpha_2+5\alpha_5=\beta_{S_4}$. On the other hand, if $\alpha_{-16}>\alpha_5$, then by Claim \ref{claim:16proofclaim7}, $\beta_S\geq 9\alpha_{-16}+2\alpha_2+8\alpha_5>8\alpha_{-16}+\alpha_2+9\alpha_5=\beta_{S_1}$. Now, if $x_2\geq3$, then by the assumption that $x_5\geq4$ and by Claim \ref{claim:16proofclaim1}, we have $\beta_S\geq5\alpha_{-16}+3\alpha_2+4\alpha_5=\beta_{S_3}$. This proves our final claim. 

Now suppose the weights are $\alpha_{-16}=\frac{1}{67}$, $\alpha_2=\frac{15}{67}$, and $\alpha_5=\frac{5}{67}$. We can check that $\kappa_0=12\alpha_{-16}+2\alpha_2+5\alpha_5=1$, achieved on the set $S_4=\{0,2,4,5,6,7,8,9,10,11,12,14,16\}$, and that this is strictly less than $\beta_{S_1}$, $\beta_{S_2}$, and $\beta_{S_3}$. By the proof of Claim \ref{claim:16proofclaim7}, any set $S$ with $x_{-16}=12,x_2=2,x_5=5$ must contain a translation of $\{0,2,4,5,6,7,8,9,10,11,12,14,16\}$, and a  so there is no possibility that a set $S$ which we did not consider, and which represents a single loop, {\em also} achieves $\beta_S=1$. This means that the walk is expected to spend an infinite amount of time in the set $\{0,2,4,5,6,7,8,9,10,11,12,14,16\}$ before ever leaving it, but there is no single loop that the walk is expected to take infinitely many times before deviating from it. 
\end{exmp}

\bibliographystyle{plain}
\bibliography{default}

\begin{thebibliography}{1}

\bibitem{Enriquez&Sabot2002}
Nathanaël Enriquez and Christophe Sabot.
\newblock Edge oriented reinforced random walks and rwre.
\newblock {\em Comptes Rendus Mathematique}, 335:941--946, 05 2002.

\bibitem{Keane&Rolles2002}
M.S. Keane and S.W.W. Rolles.
\newblock Tubular recurrence.
\newblock {\em Acta Mathematica Hungarica}, 97:207--221, 2002.

\bibitem{Sabot&Tournier2016}
Christophe Sabot and Laurent Tournier.
\newblock Random walks in dirichlet environment: an overview.
\newblock {\em Annales de la faculté des sciences de Toulouse Mathématiques},
  26:463--509, 01 2016.

\bibitem{Slonim2022Ballistic}
Daniel~J. Slonim.
\newblock Ballisticity of random walks in random environments on $\mathbb{Z}$
  with bounded jumps.
\newblock Preprint, submitted May 2022. https://arxiv.org/abs/2205.06419.

\bibitem{Slonim2021a}
Daniel~J. Slonim.
\newblock Random walks in dirichlet random environments on $\mathbb{Z}$ with
  bounded jumps.
\newblock Preprint, submitted May 2021. https://arxiv.org/abs/2104.14950.

\bibitem{Tournier2009}
Laurent Tournier.
\newblock Integrability of exit times and ballisticity for random walks in
  dirichlet environment.
\newblock {\em Electron. J. Probab.}, 14:431--451, 2009.

\end{thebibliography}

\end{document}